%!TEX TS-program =  tex
% This is the amstex file of the paper ``Discrete hypergroups associated
% with compact quantum Gelfand pairs'' by Tom H. Koornwinder, which was
% accepted for publication in the Contemporary Mathematics Volume of the
% Proceedings of the 1993 Joint Summer
% Research Conference on applications of hypergroups and related measure
% algebras, edited by W. Connett, Olivier Gebuhrer, and A. Schwartz.
%
\magnification1200
\input amstex
\documentstyle{amsppt}
\vsize 8.3truein
\NoBlackBoxes
\leftheadtext{T. H. KOORNWINDER}
\rightheadtext{COMPACT QUANTUM GELFAND PAIRS}

\define\gok{{\goth k}}
\define\CC{{\Bbb C}}
\define\RR{{\Bbb R}}
\define\ZZ{{\Bbb Z}}
\define\Zplus{\ZZ_+}
\define\al{\alpha}
\define\be{\beta}
\define\ga{\gamma}
\define\de{\delta}
\define\ep{\varepsilon}
\define\ze{\zeta}
\define\th{\theta}
\define\la{\lambda}
\define\si{\sigma}
\define\De{\Delta}
\define\FSA{{\Cal A}}
\define\FSB{{\Cal B}}
\define\FSH{{\Cal H}}
\define\FSU{{\Cal U}}
\define\FSZ{{\Cal Z}}
\define\st{such that}
\define\wrt{with respect to }
\define\id{\text{id}}
\define\supp{\text{supp}}
\define\Span{\text{Span}}
\define\ten{\otimes}
\define\rep{representation}
\define\corep{corepresentation}
\define\tr{\text{tr}\,}
\define\iy{\infty}
\define\myRe{\text{Re}\,}
\define\lan{\langle}
\define\ran{\rangle}
\define\half{\frac12}
\define\thalf{{\dsize\frac12}}
\define\wh#1{\widehat#1}
\define\wb#1{\overline{#1}}
\define\Ahat{\wh\FSA\,}
\define\ABhat{(\FSA,\FSB)\,\wh{\;}\,}
\define\AJhat{(\FSA,J)\,\wh{\;}\,}
\define\posdef{positive definite}
\define\posdefness{positive definiteness}
\define\albegade{\left(\matrix\al&\be\\ \ga&\de\endmatrix\right)}
\define\ccite#1#2{\cite{#1\rm, #2}}

\topmatter
\title
Discrete Hypergroups Associated with Compact Quantum Gelfand Pairs
\endtitle

\author
Tom H. Koornwinder
\endauthor
\address
University of Amsterdam, Faculty of Mathematics and Computer Science,
Plantage Muidergracht 24, 1018 TV Amsterdam, The Netherlands
\endaddress
\email
thk\@fwi.uva.nl
\endemail


\abstract
A discrete DJS-hypergroup is constructed in connection
with the linearization formula for
the product of two spherical elements for a quantum Gelfand pair of two
compact quantum groups. A similar construction is discussed for the case of a
generalized
quantum Gelfand pair, where the role of the quantum subgroup is taken over by
a two-sided
coideal in the dual Hopf algebra. The paper starts with a review of compact
quantum groups, with an approach in terms of so-called CQG algebras.
The paper concludes with some examples of hypergroups thus obtained.
\endabstract

\thanks
This paper appeared in
``Applications of hypergroups and related measure algebras'',
W. C. Connett, M.-O. Gebuhrer \& A. L. Schwartz (eds.),
Contemp. Math. 183, 1995, Amer. Math. Soc., pp. 213--235.
\endthanks
\endtopmatter

\document

\head 1. Introduction \endhead
Convolution algebras of $K$-biinvariant measures on a locally compact group
$G$ for a Gelfand pair $(G,K)$ were motivating examples for the introduction
of DJS-hypergroups, see for instance Jewett \cite{9}.
In the case of compact $G$, one can also associate a (discrete) dual
hypergroup with such a Gelfand pair. The convolution on this dual hypergroup
is related to the positivity of the coefficients
in the linearization formula for spherical functions.
Therefore, it is natural to expect that hypergroups may also arise in
the context of quantum groups and quantum analogues of Gelfand pairs.
In an earlier paper \cite{12} I gave many indications for this.
I introduced quantum Gelfand pairs of compact matrix quantum groups
and I showed positivity results associated with it,
but I did not go all the way to realize a hypergroup structure.

In an informal note \cite{13}, which was not widely circulated,
I showed that
there is indeed a DJS-hypergroup structure associated with these compact
quantum Gelfand pairs for the dual case.
This result is not completely trivial, because one has to replace the
involution on the Hopf $*$-algebra by another involution
(canonically determined) in order to obtain the involution for the hypergroup.

It is the purpose of the present paper to give the details of
the construction in
\cite{13}. Compared to \cite{12}, \cite{13}, the results are
formulated in terms of CQG algebras (cf.\ Koornwinder \ccite{15}{\S2},
Dijkhuizen \ccite{4}{Ch.\ 2} and, earlier with a different terminology,
Effros \& Ruan \cite{5}), which both generalize and simplify
the Woronowicz compact matrix quantum groups \cite{23}.
Furthermore, a structure of dual hypergroup is now also established
for the case of a generalized quantum Gelfand pair $(\FSA,J)$,
where $\FSA$ is a CQG algebra and $J$ is a two-sided coideal in the
algebraic linear dual $\FSA^*$ of $\FSA$.
These last results are applied to the author's case \cite{14}
where $\FSA=\FSA_q(SU(2))$ and to Noumi's case \cite{18}, where
$\FSA=\FSA_q(U(N))$ and $J$ quantizes the subgroup $SO(N)$ of $U(N)$.

This volume also contains
a very comprehensive survey by Vainerman
\cite{21} on hypergroups in relation with
quantum Gelfand pairs.
The present paper, written independently from \cite{21},
can be considered as a detailed study of a small part of the wide theory
surveyed by Vainerman.

The contents of this paper are as follows.
Section 2 gives a summary of the results on Hopf algebras and
CQG algebras needed in this paper.
Section 3 introduces a second, canonically determined
involution on a CQG algebra. The results are partly new.
In section 4 we recall from \cite{12}
the notion of a \posdef\ element in a CQG algebra.
Section 5 introduces quantum Gelfand pairs (of CQG algebras) and
generalized quantum Gelfand pairs (of a CQG algebra and a two-sided
coideal in the dual of the CQG algebra).
Part of this section recapitulates results from \cite{12}, but another
part is new.
Next, in section 6, we obtain (dual) DJS-hypergroups associated with
quantum Gelfand pairs.
Section 7 gives examples with $SU_q(2)$ and section 8 finally applies
the result to Noumi's \cite{18} generalized quantum Gelfand pair.

\remark{Acknowledgements}
I thank prof.\ H. Heyer for his pertinent question to make the allusion
to hypergroups in \cite{12} more concrete.
I also thank prof.\ G. Gasper for providing me with information about
possible analytic proofs of positivity of linearization coefficients
for certain Askey-Wilson polynomials, cf.\ \S7.
I thank Paul Floris and Erik Koelink
for carefully reading and correcting an earlier version
of this paper.
\endremark

\head 2. Preliminaries on compact quantum groups \endhead
Standard references about Hopf algebras are the books by Abe \cite{1}
and by Sweedler \cite{20},
see also the tutorial introduction to Hopf algebras in
Koornwinder \ccite{15}{\S1}.
I will summarize some definitions and properties.

\definition{Definition 2.1}
A {\sl Hopf algebra} (over $\CC$) is a complex associative unital algebra
$\FSA$ with additional linear operations
$\De\colon\FSA\to\FSA\ten\FSA$ ({\sl comultiplication}),
$\ep\colon\FSA\to\CC$ ({\sl counit}) and
$S\colon\FSA\to\FSA$ ({\sl antipode}) \st\ the following properties are
satisfied:
\roster
\item"{(a)}"
$(\De\ten\id)\circ\De=(\id\ten\De)\circ\De$\quad
({\sl coassociativity});
\item"{(b)}"
$(\ep\ten\id)\circ\De=\id=(\id\ten\ep)\circ\De$;
\item"{(c)}"
$\De$ and $\ep$ are unital algebra homomorphisms;
\item"{(d)}"
$(m\circ(S\ten\id)\circ\De)(a)=\ep(a)\,1=(m\circ(\id\ten S)\circ\De)(a)$
for all $a\in\FSA$.
\endroster
Here $m\colon\FSA\ten\FSA\to\FSA$ is the unique linear mapping \st\
$m(a\ten b)=ab$ for all $a,b\in\FSA$.
Also, $\FSA\ten\FSA$ is made into
a unital algebra \st\ $(a\ten b)\,(c\ten d)=ac\ten bd$ for all $a,b,c,d\in\FSA$.
\enddefinition

Let $\si$ be the linear endomorphism of $\FSA\ten\FSA$
\st\ $\si(a\ten b)=b\ten a$ (the {\sl flip}).
The antipode satisfies the properties $S(1)=1$, $\ep\circ S=\ep$,
$S(ab)=S(b)\,S(a)$ ($a,b\in\FSA$) and
$$
(S\ten S)\circ\De=\si\circ\De\circ S.
\tag2.1
$$

The following notation is often useful in a Hopf algebra $\FSA$.
If $a\in\FSA$ then we can choose sets of
elements $a_{(1)i}$ and $a_{(2)i}$ in $\FSA$ ($i$ running over a
finite set) \st\
$\De(a)=\sum_i a_{(1)i}\ten a_{(2)i}$.
We write this symbolically as
$$
\De(a)=\sum_{(a)}a_{(1)}\ten a_{(2)},\quad a\in\FSA.
\tag2.2
$$
Similarly, we write
$$
(\De\ten\id)(\De(a))=
\sum_{(a)}a_{(1)}\ten a_{(2)}\ten a_{(3)},\quad a\in\FSA.
$$
This notation is justified by the coassociativity property of $\De$.

A {\sl $*$-algebra} (always assumed to be unital) is a complex associative
unital algebra $\FSA$ equipped with an involutive antilinear
antimultiplicative operation $a\mapsto a^*\colon \FSA\to\FSA$.

A {\sl Hopf $*$-algebra} is a Hopf algebra which, considered as an algebra,
is a $*$-algebra \st\ the algebra homomorphisms $\De$ and $\ep$ are
homomorphisms of $*$-algebras.
Here $\FSA\ten\FSA$ is considered as a $*$-algebra \st\
$(a\ten b)^*=a^*\ten b^*$ ($a,b\in\FSA$).
In a Hopf $*$-algebra the following property is valid.
$$
S \circ * \circ S \circ *=\id.
\tag2.3
$$
By way of example let $G$ be a compact group and let $\FSA=\FSA(G)$
be the complex
linear space of all functions on $G$ which are linear combinations of
matrix elements of finite dimensional unitary matrix \rep s of $G$.
Then $\FSA$ is a commutative unital $*$-algebra under pointwise multiplication
and pointwise complex conjugation.
There is a linear embedding of $\FSA\ten\FSA$ in the space of functions on
$G\times G$ \st\ $(f\ten g)(x,y):=f(x)\,g(y)$ for $x,y\in G$ and
$f,g\in\FSA$. For $f\in\FSA$ let $\De(f)$ be the function on $G\times G$
defined by
$(\De(f))(x,y):=f(xy)$ ($x,y\in G$).
If $t_{ij}$ is a matrix element of a unitary matrix \rep\
$t=(t_{ij})_{i,j=1,\ldots,n}$
of $G$ then it follows that
$$
\De(t_{ij})=\sum_{k=1}^n t_{ik}\ten t_{kj}.
\tag2.4
$$
Hence $\De$ maps $\FSA$ into $\FSA\ten\FSA$.
Also define, for $f\in\FSA$,
that $\ep(f):=f(e)$ and $(S(f))(x):=f(x^{-1})$ ($x\in G$).
With these operations $\FSA$ becomes a commutative Hopf $*$-algebra.

Let $\FSA$ be a Hopf-algebra. Let $\FSA^*$ be its algebraic linear dual.
Then $\FSA^*$ becomes a unital algebra with identity element $\ep$
and multiplication defined by
$$
(fg)(a):=(f\ten g)(\De(a)),\quad
f,g\in\FSA^*,\;a\in\FSA.
\tag2.5
$$
Furthermore, the multiplication, unit and antipode on $\FSA$ induce
linear operations $\De\colon\FSA^*\to(\FSA\ten\FSA)^*$,
$\ep\colon\FSA^*\to\CC$  and $S\colon
\FSA^*\to\FSA^*$ as follows.
$$
\align
\De(f)(a\ten b)&:=f(ab),\quad f\in\FSA^*,\;a,b\in\FSA,\tag2.6
\\
\ep(f)&:=f(1),\quad f\in\FSA^*,\tag2.7
\\
(S(f))(a)&:=f(S(a)),\quad f\in\FSA^*,\;a\in\FSA.\tag2.8
\endalign
$$
Then $\FSA^*$, with the operations defined by formulas (2.5)--(2.8),
satisfies the axioms of a Hopf algebra, except for
slight modifications because $\De$ maps $\FSA^*$ to $(\FSA\ten\FSA)^*$
rather than $\FSA^*\ten\FSA^*$.
If $\FSA$ is moreover a Hopf $*$-algebra then $\FSA^*$ becomes a $*$-algebra
(almost a Hopf $*$-algebra) with involution given by
$$
f^*(a):=\wb{f((S(a))^*)},\qquad f\in\FSA^*,\;a\in\FSA.
\tag2.9
$$

Consider the example $\FSA=\FSA(G)$ ($G$ a compact group).
Then a complex regular Borel measure $\mu$ on
$G$ determines a linear functional
$\mu\colon f\mapsto \int_G f(x)\,d\mu(x)$ on $\FSA$ and $\mu$
is the unique regular measure which gives rise to this linear functional.
Then the product of two regular Borel measures $\mu,\nu$ on $G$
considered as linear
functionals on $\FSA$ is the linear functional on $\FSA$ corresponding
to the convolution product
$\mu*\nu$, while $\mu^*(f)=\wb{\int_G \wb{f(x^{-1})}\,d\mu(x)}$.

Let $\FSA$ be a Hopf algebra.
We can define left and right actions of the unital algebra $\FSA^*$
on $\FSA$ as follows. 
$$
\align
&f.a:=(\id\ten f)(\De(a))=
\sum_{(a)}f(a_{(2)})\,a_{(1)},\quad f\in\FSA^*,\;a\in\FSA,\tag2.10
\\
&a.f:=(f\ten\id)(\De(a))=
\sum_{(a)}f(a_{(1)})\,a_{(2)},\quad f\in\FSA^*,\;a\in\FSA.\tag2.11
\endalign
$$
Then 
$$
(fg).a=f.(g.a),\quad a.(fg)=(a.f).g\qquad(f,g\in\FSA^*,\;a\in\FSA).
$$
Also
$$
f(g.a)=(fg)(a)=g(a.f)\qquad
(f,g\in\FSA^*,\;a\in\FSA).
\tag2.12
$$
In particular,
$$
\ep(f.a)=f(a)=\ep(a.f)\quad(f\in\FSA^*,\;a\in\FSA).
\tag2.13
$$
If $f\in\FSA^*$ \st\ $\De(f)=\sum_{(f)}f_{(1)}\ten f_{(2)}\in\FSA^*\ten\FSA^*$
then
$$
f.(ab)=\sum_{(f)} (f_{(1)}.a)\,(f_{(2)}.b),\quad
(ab).f=\sum_{(f)} (a.f_{(1)})\,(b.f_{(2)})\qquad (a,b\in\FSA).
\tag2.14
$$
A {\sl matrix \corep} of a Hopf algebra $\FSA$
is a matrix $t=(t_{ij})_{i,j=1,\ldots,n}$
with entries in $\FSA$ \st\ (2.4) is valid and also
$$
\ep(t_{ij})=\de_{ij}.
\tag2.15
$$
A matrix \corep\ $t$ of a Hopf $*$-algebra is called {\sl unitary} if
$$
S(t_{ij})=(t_{ji})^*.
\tag2.16
$$
In the example of $\FSA=\FSA(G)$, where $G$ is a compact group,
a square matrix $t$ with elements in $\FSA$ is a unitary matrix \rep\ of $G$
if and only if it is a unitary matrix \corep\ of $\FSA$.

Two matrix \corep s $s$ and $t$ of a Hopf algebra,
both of the same dimension $n$, are called {\sl equivalent} if there is a
complex invertible $n\times n$ matrix $B$ \st\
$B\,s=t\,B$. 

A matrix \corep\ $t$ (of dimension $n$) of a Hopf algebra is called
{\sl irreducible} if $t$ is not equivalent to a matrix \corep\ of the form
$\left(\matrix*&*\\0&*\endmatrix\right)$.

If $t=(t_{ij})$ is a matrix \corep\ of $\FSA$ then so is its {\sl contragredient
\corep} $t'$, where $(t')_{ij}:=S(t_{ji})$. Then the \corep\ $t'':=(t')'$
has matrix elements $(t'')_{ij}=S^2(t_{ij})$.

In Koornwinder \ccite{15}{\S2} the following two definitions were given.
A {\sl CQG algebra} is a Hopf $*$-algebra which is the linear span of the
matrix elements of its unitary matrix \corep s.
A {\sl CMQG algebra} is a Hopf $*$-algebra $\FSA$ which, as a unital algebra,
is generated by the matrix elements of a certain unitary matrix \corep\ of
$\FSA$.
It follows easily (cf.\ \cite{15}) that a Hopf $*$-algebra is
a CMQG algebra if and only if it is a CQG algebra which, as an algebra,
is finitely generated.
Clearly, the Hopf $*$-algebra $\FSA(G)$ ($G$ a compact group) is a CQG algebra.
It is a CMQG algebra if and only if $G$ is isomorphic to a closed
subgroup of a unitary group $U(n)$, which, in its turn, is equivalent
to the fact that $G$ is a compact Lie group.

The theory of CQG algebras was developed by the author
\ccite{15}{\S2} jointly with M.~S. Dijkhuizen
\cite{4}. It is intended as an alternative approach to the compact
matrix quantum groups introduced by Woronowicz
\cite{23}, \cite{25}, \cite{26}.
Prior to \cite{15}, CQG algebras were introduced under a different name
by Effros \& Ruan \cite{5}. However, their further development of the
theory is very different from the approach in \cite{15}.

Let $\FSA$ be a CQG algebra. We quote some further results from
\ccite{15}{\S2}.
Denote by $\Ahat$ the collection of all equivalence classes of
irreducible unitary matrix \corep s of $\FSA$. Choose,
for each $\al\in\Ahat$, a unitary matrix \corep\
$(t_{ij}^\al)_{i,j=1,\ldots,d_\al}$ belonging to class $\al$.
Then the set of all $t_{ij}^\al$ forms a
basis of $\FSA$.
Let  $e$ be the element of $\Ahat$ for which $t^e$ is the one-dimensional
matrix \corep\ $(1)$ ({\sl trivial \corep}).

Define the {\sl Haar functional} $h$ on $\FSA$ as the linear mapping
$h\colon\FSA\to\CC$ \st
$$
h(t_{ij}^\al)=0\quad(\al\ne e),\qquad
h(1)=1.
$$
Then it is immediately verified that
$$
(h\ten\id)(\De(a))=h(a)\,1=(\id\ten h)(\De(a)),\quad a\in\FSA,
$$
while it is a deeper result that
$$
h(aa^*)>0\quad\hbox{if $a\ne0$.}
$$
In the example $\FSA=\FSA(G)$ ($G$ a compact group) we have
$h(f)=\int_G f(x)\,dx$, where $dx$ is the normalized Haar measure on $G$.

\proclaim{Proposition 2.2 {\rm(cf.\ \ccite{15}{\S2})}}
Let $r$ be an irreducible unitary matrix \corep\ of a CQG algebra $\FSA$.
Then the \corep s $r$ and $r''$ are equivalent. Moreover, if $F$
is a complex invertible matrix \st\ $F\,r=r''\,F$ then
$\tr F\ne0\ne\tr F^{-1}$ and $F$ can be uniquely chosen \st\
$\tr F=\tr F^{-1}>0$. Then $F$ is a \posdef\ matrix.
\endproclaim

If $\al\in\Ahat$ then write $F_\al$ for the normalized operator $F$
associated with $t^\al$ by Proposition 2.2.
The {\sl quantum Schur orthogonality relations} read as follows.
$$
 h(t_{kl}^\be\,(t_{ji}^\al)^*)=
\de_{\al\be}\,\de_{kj}\,{(F_\al)_{il}\over\tr F_\al}\,,\quad
h((t_{lk}^\be)^*\,t_{ij}^\al)=
\de_{\al\be}\,\de_{kj}\,{(F_\al^{-1})_{il}\over\tr F_\al^{-1}}\qquad
(\al,\be\in\Ahat).
\tag2.17
$$
All these results are essentially due to Woronowicz \cite{23},
but they are proved somewhat differently in the CQG algebra approach
of \ccite{15}{\S2}.

I will now point out how CQG algebras are related to the compact matrix quantum
groups of Woronowicz, cf.\ \ccite{15}{\S2.3, \S2.5}.
For each element $a$ of a CQG algebra $\FSA$ put
$\|a\|:=\sup_{\pi} \|\pi(a)\|$.
Here $\pi$ runs through all $*$-\rep s on Hilbert spaces
of the $*$-algebra $\FSA$. Then, by making essential use of the fact that $\FSA$
is spanned by matrix elements of unitary \corep s, it can be
shown first that $\|a\|<\iy$ for all $a\in\FSA$, and next that the seminorm
$\|\,.\,\|$ is in fact a norm on $\FSA$. This norm satisfies
$\|aa^*\|=\|a\|^2$ ($C^*$-norm). Denote the completion of $\FSA$ \wrt
this norm by $A$. Then $A$ is a unital $C^*$-algebra.
Equip the algebraic tensor product $A\ten A$ with the $C^*$-norm
$\|a\|:=\sup_{\pi_1,\pi_2}\|(\pi_1\ten\pi_2)(a)\|$, where
$a\in A\ten A$ and $\pi_1, \pi_2$ run through all $*$-\rep s on Hilbert spaces
of $A$. Then $\De\colon\FSA\to\FSA\ten\FSA$ continuously extends to 
a $C^*$-homomorphism from $A$ to the completion of $A\ten A$ \wrt this norm.

It can be deduced from Woronowicz \cite{23} that the
{\sl compact matrix quantum groups} defined there,
can be characterized as pairs $(A,\FSA)$,
where $A$ is a $C^*$-algebra, $\FSA$ is a CMQG algebra and
a dense $*$-subalgebra of $A$, and
$\De\colon\FSA\to\FSA\ten\FSA$ continuously extends to a
$C^*$-homomorphism from $A$ to a suitable $C^*$-completion of $A\ten A$.
Thus a CQG algebra $\FSA$ can be brought in correspondence with a
Woronowicz compact matrix quantum group $(A,\FSA)$
if and only if $\FSA$ is a CMQG algebra,
but $\FSA$ does not determine $A$ uniquely.

\head 3. A second involution for CQG algebras \endhead
Let $\FSA$ be a Hopf algebra.
A nonzero linear functional $f$ on $\FSA$ is called {\sl multiplicative}
if $f(ab)=f(a)\,f(b)$ for all $a,b\in\FSA$. Then also $f(1)=1$.
A multiplicative linear functional on $\FSA$ can equivalently be
characterized as a nonzero element $f\in\FSA^*$ \st\
$$
\De(f)=f\ten f.
\tag3.1
$$
Such elements are called {\sl group-like}. In the example $\FSA=\FSA(G)$
the point evaluations $f\mapsto f(x)\colon \FSA\to\CC$ ($x\in G$)
yield group-like elements of $\FSA^*$.

If $f\in\FSA^*$ is group-like then
$$
f\,S(f)=\ep=S(f)\,f.
\tag3.2
$$
Indeed,
$$
\multline
(f\,S(f))(a)=(f\ten f)((\id\ten S)(\De(a)))=
(\De(f))((\id\ten S)(\De(a)))
\\
=f((m\circ(\id\ten S)\circ\De)(a))=f(\ep(a)\,1)=\ep(a)\,f(1)=\ep(a).
\endmultline
$$
and similarly for the other identity. Here we used (2.5), (2.8),
(3.1), (2.6) and Definition 2.1(d).
If $f$ is group-like then so is $Sf$.

If $f\in\FSA^*$ is group-like then let the linear mapping
$\tau_f\colon\FSA\to\FSA$ be defined by
$$
\tau_f(a):=S(f).a.f=
\sum_{(a)} f(a_{(1)})\,a_{(2)}\,(S(f))(a_{(3)}),\qquad a\in\FSA.
\tag3.3
$$
We call the mappings $\tau_f$ {\sl inner automorphisms} of $\FSA$,
since the next Proposition shows that $\tau_f$ is a Hopf algebra
automorphism and since, if $\FSA=\FSA(G)$ and $\chi(f):=f(x)$
for some $x\in G$, then $(\tau_\chi(f))(y)=f(xyx^{-1})$\quad($y\in G$).

\proclaim{Proposition 3.1}
Let $\FSA$ be a Hopf algebra and let $f\in\FSA^*$ be group-like.
Then
$$
\align
&\tau_{S(f)}\circ\tau_f=\id=\tau_f\circ\tau_{S(f)},\tag3.4
\\
&f.(ab)=(f.a)\,(f.b),\quad
(ab).f=(a.f)\,(b.f)\qquad (a,b\in\FSA),\tag3.5
\\
&\tau_f(ab)=\tau_f(a)\,\tau_f(b)\quad (a,b\in\FSA),\tag3.6
\cr
&\De\circ\tau_f=(\tau_f\ten\tau_f)\circ\De.\tag3.7
\endalign
$$
\endproclaim

\demo{Proof}
Formula (3.4) follows immediately from (3.3) and (3.2).
Formula (3.5) follows from (2.14).
Then (3.5) immediately yields (3.6).
Finally, (3.7) is proved as follows.
$$
\align
\De(\tau_f(a))&=
\De\Bigl(\sum_{(a)}f(a_{(1)})\,a_{(2)}\,(S(f))(a_{(3)})\Bigr)
\\
&=\sum_{(a)}f(a_{(1)})\,a_{(2)}\,a_{(3)}\,(S(f))(a_{(4)})
\\
&=\sum_{(a)}f(a_{(1)})\,a_{(2)}\,\ep(a_{(3)})\,a_{(4)}\,(S(f))(a_{(5)})
\\
&=\sum_{(a)}f(a_{(1)})\,a_{(2)}\,(S(f))(a_{(3)})\,f(a_{(4)})\,
a_{(5)}\,(S(f))(a_{(6)})
\\
&=\sum_{(a)}\tau_f(a_{(1)})\,\tau_f(a_{(2)}).
\endalign
$$
Here we used (3.3), Definition 2.1(c), (3.2) and
(2.5).
\qed\enddemo

\remark{Remark 3.2}
Let $f\in\FSA^*$ be group-like.
If $t=(t_{ij})$ is a matrix \corep\ of the Hopf algebra $\FSA$ then put
$r_{ij}:=\tau_f(t_{ij})$ and $r:=(r_{ij})$.
Then $r$ is a matrix \corep\ of $\FSA$
which is equivalent to $t$. Indeed,
$$
r_{ij}=S(f).t_{ij}.f=\sum_{k,l}(f)(t_{ik})\,t_{kl}\,(S(f))(t_{lj})
$$
by (3.3) and (2.4).
Furthermore, the matrices $\bigl(f(t_{ij})\bigr)$ and
$\bigl((S(f))(t_{ij})\bigr)$ are inverse to each other by
(3.2) and (2.4).
This settles the equivalence.
\endremark\medpagebreak

Let now $\FSA$ be a Hopf $*$-algebra and let $f\in\FSA^*$ be group-like.
Define the mapping $a\mapsto a^-\colon\FSA\to\FSA$ by
$$
a^-:=(\tau_{S(f)}\circ*\circ\tau_f)(a).
\tag3.8
$$
Then it follows immediately that $a\mapsto a^-$ is an involutive anti-linear
mapping and that, by Proposition 3.1, the Hopf algebra $\FSA$,
together with the involution $a\mapsto a^-$, has the structure of
a Hopf $*$-algebra.

Let $\FSA$ be a CQG algebra.
The intertwining operators $F_\al$ ($\al\in\Ahat$)
give rise to a remarkable family of multiplicative
linear functionals $f_z$ ($z\in\CC$)
on $\FSA$ which were introduced by Woronowicz \cite{23},
see also \ccite{15}{\S2.4}.
We will summarize their properties.
Define for each $z\in\CC$ the linear functional $f_z$ on $\FSA$ \st
$$
f_z(t_{ij}^\al):=(F_\al^z)_{ij},\quad \al\in\Ahat,\;i,j=1,\ldots,d_\al.
\tag3.9
$$
Here arbitrary complex powers $F_\al^z$ of the \posdef\ matrix $F_\al$
are defined in an evident way.
Note that the definition of $f_z$ is independent of the choice
of the unitary matrix \corep\ representing $\al\in\Ahat$.

\proclaim{Proposition 3.3}
\roster
\item"{(a)}"
$f_z\,f_{z'}=f_{z+z'}$,\quad
$f_0=\ep$.
\item"{(b)}"
$f_z(a)\,f_z(b)=f_z(ab)$,\quad
$f_z(1)=1$.
\item"{(c)}"
$S(f_z)=f_{-z}$,\quad
$(f_z)^*=f_{\wb z}$.
\item"{(d)}"
$S^2(a)=f_{-1}.a.f_1$.
\item"{(e)}"
$h(ab)=h(b\,(f_1.a.f_1))$.
\item"{(f)}"
For each $a\in\FSA$ the function $z\mapsto f_z(a)$ is an entire analytic
function and there are constants $M>0$ and $\mu\in\RR$ \st\
$|f_z(a)|\le M\,e^{\mu\,\myRe z}$.
\endroster
\endproclaim

Thus, for each $z\in\CC$, we have a group-like element $f_z\in\FSA$
which gives rise to an inner automorphism
$\tau_{f_z}\colon a\mapsto f_{-z}.a.f_z$ of $\FSA$.
We will consider the corresponding involution defined by (3.8).
First we need a few simple properties concerning the action of
$\FSA^*$ on a Hopf $*$-algebra $\FSA$.

\proclaim{Lemma 3.4}
Let $\FSA$ be a Hopf $*$-algebra. Let $f\in\FSA^*$, $a,b\in\FSA$. Then
\roster
\item"{(a)}"
$f(a^*)=\wb{(S(f))^*(a)}$.
\item"{(b)}"
$S(S(f).a)=S(a).f$,\quad
$S(a.S(f))=f.S(a)$.
\item"{(c)}"
$(f.a)^*=(S(f))^*.a^*$,\quad
$(a.f)^*=a^*.(S(f))^*$.
\endroster
\endproclaim

\demo{Proof}
For the proof of (a) note that
$$
\wb{(S(f))^*(a)}=
(S(f))((S(a))^*)=
f((S\circ*\circ S)(a))=
f(a^*),
$$
where we used (2.9), (2.8) and (2.3).

The first identity in (b) is proved as follows.
$$
\split
S(S(f).a)=S(\id\ten S(f))\,(\De(a))&
=(\id\ten f)(S\ten S)(\De(a))
\\
&\qquad
=(f\ten\id)\,\De(S(a))=
S(a).f.
\endsplit
$$
Here we used (2.10), (2.8), (2.1) and (2.11).
The second identity can be proved in a similar way.

For the proof of (c) we write:
$$
(f.a^*)^*=
\sum_{(a)}\wb{f(a_{(2)}^*)}\,a_{(1)}=
\sum_{(a)}(S(f))^*(a_{(2)})\,a_{(1)}=
(S(f))^*.a.
$$
Here we used (2.10) and part (a) of the Lemma.
This settles the first identity in (c). The second identity is proved in a
similar way.
\qed\enddemo\medpagebreak

Let $\FSA$ be a CQG algebra,
let $z\in\CC$ and let the involution $a\mapsto a^-$ be defined by
(3.8) with $f:=f_z$. Then, by Lemma 3.4(c),
$$
f_z.(f_{-z}.a.f_z)^*.f_{-z}=
f_z.f_{\wb z}.a^*.f_{-\wb z}.f_{-z}=f_{2\myRe z}.a^*.f_{-2\myRe z}\qquad
(a\in\FSA).
$$
Hence
$$
a^-=f_{2\myRe z}.a^*.f_{-2\myRe z}=\tau_{f_{-2\myRe z}}(a^*).
\tag3.10
$$
Let $t=(t_{ij})$ be a unitary matrix \corep\ of $\FSA$ and put
$(t^*)_{ij}:=(t_{ij})^*$ and $(t^-)_{ij}:=(t_{ij})^-$.
Then $t^*$ and $t^-$ are equivalent matrix \corep s of $\FSA$ because of
(3.10) and Remark 3.2.
In general, the \corep\ $t^*$ is not unitary.
We want to determine $z$ \st\ $t^-$ is unitary,
independent of the choice of $t$.
So we want that
$$
S(t_{ij}^-)=(t_{ji}^-)^*.
\tag3.11
$$
Now, on the one hand we have
$$
\multline
S(t_{ij}^-)=
S(f_{2\myRe z}.t_{ij}^*.f_{-2\myRe z})=
S(S(f_{-2\myRe z}).S(t_{ji}).S(f_{2\myRe z}))
\\
=f_{2\myRe z}.S^2(t_{ji}).f_{-2\myRe z}=
f_{2\myRe z}.f_{-1}.t_{ji}.f_1.f_{-2\myRe z}=
f_{-1+2\myRe z}.t_{ji}.f_{1-2\myRe z}.
\endmultline
$$
Here we used (3.10), (2.16), Proposition 3.3(c),
Lemma 3.4(b) and Proposition 3.3(d) and (a).
On the other hand,
$$
(t_{ji}^-)^*=(f_{2\myRe z}.t_{ji}^*.f_{-2\myRe z})^*=
((S(f_{2\myRe z}))^*.t_{ji}.(S(f_{-2\myRe z}))^*=
f_{-2\myRe z}.t_{ji}.f_{2\myRe z}.
$$
Here we used (3.10), Lemma 3.4(c) and
Proposition 3.3(c).
Thus (3.11) is satisfied if $\myRe z=\frac14$.
Because of (3.10), we may as well take $z$ real and equal to~$\frac14$.

We summarize the results in the following theorem.

\proclaim{Theorem 3.5}
Let $\FSA$ be a CQG algebra. Then the mapping $a\mapsto a^-\colon\FSA\to\FSA$
defined by
$$
a^-:=f_\half.a^*.f_{-\half}\qquad(a\in\FSA)
\tag3.12
$$
is an involution on $\FSA$ \st
$$
(f_{-\frac14}.a.f_{\frac14})^*=f_{-\frac14}.a^-.f_{\frac14}.
$$
Then $\FSA$ as a Hopf algebra,
together with this involution, becomes a Hopf $*$-algebra.
Furthermore, if $\FSA$ is considered as Hopf $*$-algebra \wrt the involution
$a\mapsto a^*$ and if $t$ is a unitary matrix \corep\ of $\FSA$ then
the matrix \corep\ $t^-$ of $\FSA$, equivalent to
the matrix \corep\ $t^*$, is unitary.
\endproclaim

We can now use (2.9) in order to define an involution $f\mapsto f^-$
on $\FSA^*$ which corresponds to the involution $a\mapsto a^-$ on $\FSA$:
$$
f^-(a):=\wb{f( (S(a))^- )},\qquad f\in\FSA^*,\;a\in\FSA.
\tag3.13
$$
Then
$$
f^-=f_\half\,f^*\,f_{-\half},\qquad f\in\FSA^*.
\tag3.14
$$
For the proof note that, with $a\in\FSA$,
$$
\multline
f^-(a)=\wb{f((S(a) )^- )}=
\wb{f(f_\half.(S(a) )^* .f_{-\half} ) }=
\wb{(f_{-\half}\,f\,f_\half )((S(a) )^* ) }
\\
=(f_{-\half}\,f\,f_\half)^*(a)=
(f_\half\,f^*\,f_{-\half})(a).
\endmultline
$$
Here we used (3.13), (3.12), (2.12), (2.9) and
Proposition 3.3(c).

The following version of the quantum Schur orthogonality
relations (2.17) involves the involution $a\mapsto a^-$.
$$
h\bigl(t_{ij}^\al\,(t_{kl}^\be)^-\bigr)=
{\de_{\al\be}\,f_\half(t_{ik}^\al)\,f_\half(t_{lj}^\al)\over\tr F_\al},\quad
h\bigl((t_{kl}^\be)^-\,t_{ij}^\al\bigr)=
{\de_{\al\be}\,f_{-\half}(t_{ik}^\al)\,f_{-\half}(t_{lj}^\al)\over\tr F_{-\al}}
\,.
\tag3.15
$$
Let us prove the first equality. The proof of the second one is analogous.
The case $\al\ne\be$ is clear from (2.17).
For $\al=\be$ we have
$$
h\bigl(t_{ij}^\al\,(t_{kl}^\al)^-\bigr)=
h\bigl(t_{ij}^\al\,f_\half.(t_{kl}^\al)^*.f_{-\half}\bigr)=
\sum_{r,s}f_{-\half}((t_{kr}^\al)^*)\,f_\half((t_{sl}^\al)^*)\,
h\bigl(t_{ij}^\al\,(t_{rs}^\al)^*\bigr),
$$
where we used (3.12), (2.10), (2.11) and (2.4).
Hence
$$
\multline
(\tr F_\al)\,h\bigl(t_{ij}^\al\,(t_{kl}^\al)^-\bigr)
=\sum_{r,s} f_\half(t_{rk}^\al)\,f_{-\half}(t_{ls}^\al)\,\de_{ir}\,(F_\al)_{sj}
\\
=\sum_s f_\half(t_{ik}^\al)\,f_{-\half}(t_{ls}^\al)\,f_1(t_{sj}^\al)
=f_\half(t_{ik}^\al)\,f_\half(t_{lj}^\al),
\endmultline
$$
where we used (2.16), (2.8), Proposition 3.3(c), (2.17),
(3.9), (2.4), (2.5) and Proposition 3.3(a).
As a special case of (3.15) we have
$$
h\bigl(t_{ij}^\al\,(t_{ij}^\al)^-\bigr)=
{f_\half(t_{ii}^\al)\,f_\half(t_{jj}^\al)\over\tr F_\al}>0,
\tag3.16
$$
where the last inequality follows because $\tr F_\al>0$ and
$f_\half(t_{ii}^\al)=(F_\al^\half)_{ii}>0$ by \posdefness\ of $F_\al$.

\head 4. Positive definite elements \endhead
The notion of a \posdef\ element was introduced in \ccite{12}{\S7}
for compact matrix quantum groups. The properties proved there remain true
for \posdef\ elements of CQG algebras and will be recapitulated below.
Some further properties can be formulated with the aid of the involution
$a\mapsto a^-$ defined in (3.12).

Let $\FSA$ be a CQG algebra.
An element $a\in\FSA$ is called {\sl \posdef} if
$$
(f^*\ten f)(\De(a))\ge 0\quad
\hbox{for all $f\in\FSA^*$.}
$$
In the example $\FSA=\FSA(G)$ the notion of \posdef\
element of $\FSA$ coincides with the notion
of \posdef\ function on $G$ (as can be seen from the next proposition).
The following properties were proved in \ccite{12}{\S7}.

\proclaim{Proposition 4.1}
Let $\FSA$ be a CQG algebra.
\roster
\item"{(a)}"
Let $(t_{ij})_{i,j=1,\ldots,n}$ be a unitary matrix \corep\ of $\FSA$.
Let also $(a_{ij})_{i,j=1,\ldots,n}$ be a \posdef\ complex hermitian matrix.
Then the element $\sum_{i,j=1}^n a_{ij}\,t_{ij}$ is \posdef.
\item"{(b)}"
Write $a\in\FSA$ as a linear combination of the $t_{ij}^\al$:
$$
a=\sum_{\al\in\Ahat}\;\sum_{i,j=1}^{d_\al} a_{ij}^\al\,t_{ij}^\al.
\tag4.1
$$
Then $a$ is \posdef\
if and only if,
for each $\al\in\Ahat$, the matrix $(a_{ij}^\al)$ is positive semi-definite.
\item"{(c)}"
If $a\in\FSA$ is \posdef\ then $\ep(a)\ge0$ and $S(a)=a^*$.
\item"{(d)}"
The element $1$ of $\FSA$ is \posdef.
\item"{(e)}"
If $a$ and $b$ are \posdef\ elements of $\FSA$ then $ab$ is \posdef.
\endroster
\endproclaim

Since $S^2\ne\id$ in general, \posdefness\ of $a$
will not imply that $a^*$ is \posdef.
The involution $a\mapsto a^-$ introduced in (3.12) has more pleasant
properties \wrt \posdefness.

\proclaim{Proposition 4.2}
Let $\FSA$ be a CQG algebra.
Then the following properties hold.
\roster
\item"{(a)}"
If $a$ is \posdef\ then so is $a^-$.
\item"{(b)}"
If both $a$ and $a^*$ are \posdef\ then $a^-=a^*$.
\endroster
\endproclaim

\demo{Proof}
For the proof of (a) expand $a$ as in (4.1).
Then all matrices $(a_{ij}^\al)$ are positive semi-definite.
We obtain from (4.1) that
$$
a^-=\sum_{\al\in\Ahat}\;\sum_{i,j=1}^{d_\al}
\wb{a_{ij}^\al}\,(t_{ij}^\al)^-.
$$
Then the matrix \corep s $((t_{ij}^\al)^-)$ are unitary by Theorem 3.5
and the matrices $(\,\wb{a_{ij}^\al}\,)$ are positive
semi-definite.
Hence, by Proposition 4.1(a), the element $a^-$ is \posdef.

Next we prove (b). Since $a$ and $a^*$ are \posdef, it follows from
Proposition 4.1(c) that $S(a)=a^*$ and $S(a^*)=a$.
Hence $S^2(a^*)=a^*$, so $f_{-1}.a^*.f_1=a^*$ by Proposition 3.3(d).
Thus, by iteration, $a^*.f_z=f_z.a^*$ for $z=1,2,\ldots\;$.
Hence, by (2.12), we have for all $z=1,2,\ldots\;$
that
$$
f_z(g.a^*)=f_z(a^*.g),\quad g\in\FSA^*.
\tag4.2
$$
Then it follows by Proposition 3.3(f) and by some function theoretic
argument (for instance using Carlson's theorem, cf.\ \ccite{15}{\S2.4}) that
(4.2) is valid for all $z\in\CC$.
Once more by (2.12), we obtain that $a^*.f_z=f_z.a^*$ for all $z\in\CC$.
Hence, $f_\half.a^*.f_{-\half}=a^*$, which proves (b).
\qed\enddemo

\head 5. Quantum Gelfand pairs \endhead
In the author's paper \ccite{12}{\S7} quantum Gelfand pairs of compact matrix
quantum groups were introduced.
Here these results will be briefly recapitulated, but reformulated now
in terms of CQG algebras.
We will also consider a generalized notion of quantum Gelfand pair, where
we deal with a CQG algebra $\FSA$ and a two-sided coideal $J$ in $\FSA^*$.
In Vainerman \cite{21} a different definition of quantum Gelfand pair
was given, but it is equivalent to the definition below.

Let $\FSA$ be a CQG algebra.
A {\sl sub-CQG algebra} of $\FSA$ is a pair $(\FSB,\Psi)$ with $\FSB$ a
CQG algebra and $\Psi\colon\FSA\to\FSB$ a surjective homomorphism
of Hopf $*$-algebras.
For example,
let $G$ be a compact group with closed subgroup $H$,
write $\FSA:=\FSA(G)$, $\FSB:=\FSA(H)$
and put $\Psi\colon f\mapsto f|_H\colon \FSA\to\FSB$
(restriction of functions on $G$ to functions on $H$). Then
$(\FSB,\Psi)$ is a sub-CQG algebra of $\FSA$.

Let $(\FSB,\Psi)$ be a sub-CQG algebra of a CQG algebra $\FSA$.
An element $a\in\FSA$ is called {\sl biinvariant} \wrt $(\FSB,\Psi)$ if
$$
(\Psi\ten\id)(\De(a))=1_\FSB\ten a\quad\text{ and}\quad
(\id\ten\Psi)(\De(a))=a\ten 1_\FSB.
$$
These biinvariant elements form a unital $*$-subalgebra of $\FSA$.

A pair of a CQG algebra $\FSA$ and a sub-CQG algebra $(\FSB,\Psi)$ of $\FSA$
is called a {\sl quantum Gelfand pair of CQG algebras}
if, for each $\al\in\Ahat$, the elements in
$\Span\{t_{ij}^\al\}$ which are biinvariant \wrt $(\FSB,\Psi)$,
form a subspace of dimension 0 or 1. From
now on we will suppress the surjective Hopf $*$-algebra homomorphism
in our notation: we will speak about the quantum Gelfand pair
$(\FSA,\FSB)$ and about $\FSB$-biinvariant elements in $\FSA$.
If $(\FSA,\FSB)$ is a quantum Gelfand pair then the set of all $\al\in\Ahat$
for which the $\FSB$-biinvariant elements in
$\Span\{t_{ij}^\al\}$ form a subspace of dimension 1, will be
denoted by $(\FSA,\FSB)\,\wh{\;}$. We will assume that, for each
$\al\in\ABhat$, the unitary matrix \corep\ $(t_{ij}^\al)$ is chosen \st\
$t_{11}^\al$ is $\FSB$-biinvariant (this is always possible).
This element is called a {\sl spherical element} for the quantum Gelfand pair
$(\FSA,\FSB)$.
Now the following proposition (cf.\ \ccite{12}{\S7})
follows immediately from Proposition 4.1.

\proclaim{Proposition 5.1}
Let $(\FSA,\FSB)$ be a quantum Gelfand pair of CQG algebras.
\roster
\item"{(a)}"
If $\al\in\ABhat$ then $t_{11}^\al$ is the unique element $a$ of
$\Span\{t_{ij}^\al\}$ which is $\FSB$-biinvariant and which satisfies
$\ep(a)=1$.
\item"{(b)}"
If $\al\in\ABhat$ then $t_{11}^\al$ is \posdef.
\item"{(c)}"
The element $e$ of $\Ahat$ belongs to $\ABhat$.
\item"{(d)}"
Let $\FSZ$ be the unital $*$-subalgebra of $\FSB$-biinvariant elements
of $\FSA$. Then
$$
\FSZ=\sum_{\al\in\ABhat}\CC\,t_{11}^\al.
$$
\item"{(e)}"
We have
$$
t_{11}^\al\,t_{11}^\be=\sum_{\ga\in\ABhat} c_{\al\be}(\ga)\,t_{11}^\ga,\quad
\al,\be\in\ABhat,
\tag5.1
$$
where only finitely many terms in the sum are nonzero and
$c_{\al\be}(\ga)\ge0$ for all $\al,\be,\ga\in\ABhat$.
\endroster
\endproclaim

Let $\FSA$ be a CQG algebra.
If $\al\in\Ahat$ then write $\wb\al$ for 
the element of $\Ahat$ for
which the irreducible unitary matrix \corep $\bigl((t_{ij}^\al)^-\bigr)$ is a
representative. Recall (cf.\ Theorem 3.5) that this last matrix \corep\
is equivalent to the matrix \corep\  $\bigl((t_{ij}^\al)^*\bigr)$.
Note that the \corep s $\bigl((t_{ij}^\al)^-\bigr)$
and $(t_{ij}^{\wb\al})$ are equivalent  but not necessary equal.
Observe that $\wb{\wb\al}=\al$\quad($\al\in\Ahat$).

The next Proposition describes how $\al\mapsto\wb\al$ acts on $\ABhat$.

\proclaim{Proposition 5.2}
Let $(\FSA,\FSB)$ be a quantum Gelfand pair of CQG algebras.
Let $\al\in\ABhat$. Then $\wb\al\in\ABhat$,
the element $(t_{11}^\al)^*$ is spherical and
$$
(t_{11}^\al)^-=(t_{11}^\al)^*=t_{11}^{\wb\al}.
\tag5.2
$$
\endproclaim

\demo{Proof}
The element $(t_{11}^\al)^*$ is $\FSB$-biinvariant and satisfies
$\ep((t_{11}^\al)^*)=1$, because $t_{11}^\al$ has such properties.
Hence $\wb\al\in\ABhat$ and the last equality in (5.2)
is satisfied.
Then, by Proposition 4.2(b), the first equality in (5.2)
also holds, because both $t_{11}^\al$ and $(t_{11}^\al)^*$ are
\posdef.
\qed\enddemo\medpagebreak

The reader may now continue in the next section, where the hypergroup
structure following from Propositions 5.1 and 5.2 is described.
Below we introduce a generalized notion of quantum Gelfand pair.

Let $\FSA$ be a CQG algebra.
Let $J$ be a {\sl two-sided coideal} of $\FSA^*$, i.e., a linear subspace \st\
$\De(J)\i J\ten \FSA^*+\FSA^*\ten J$ and $f(1)=0$ for all $f\in J$.
An element $a\in\FSA$ is called {\sl $J$-biinvariant} if
$f.a=0=a.f$ for all $f\in J$.
The $J$-biinvariant elements form a unital subalgebra of $\FSA^*$.
This follows from (2.14).

If $t=(t_{ij})$ is a unitary matrix \corep\ of a Hopf $*$-algebra $\FSA$
and if we put
$\pi_{ij}(f):=f(t_{ij})$ for $f\in\FSA^*$, then
$\pi\colon f\mapsto(\pi_{ij}(f))$ is a matrix $*$-\rep\ of
the $*$-algebra $\FSA^*$. If $\FSA$ is a CQG algebra and $\al\in\Ahat$
then we write $\pi^\al$ for the $*$-\rep\ of $\FSA^*$  thus corresponding
to the unitary \corep\ $t^\al$ of $\FSA$. Let $\FSH^\al$ denote the
\rep\ space on which the \rep\ $\pi^\al$ is acting.
It can be identified with $\CC^{d_\al}$.
Let $e_1,e_2,\ldots,e_{d_\al}$ be the standard basis of $\FSH^\al$.

Let $\FSA$ be a CQG algebra and let $J$ be a two-sided coideal of $\FSA^*$
\st\ $J=J^*$. We call the pair $(\FSA,J)$ a {\sl generalized
quantum Gelfand pair} if, for each $\al\in\Ahat$, the
$J$-invariant elements in $\FSH^\al$ form a subspace of dimension 0 or 1.
Let $\AJhat$ denote the set of all $\al\in\Ahat$ for which this dimension is 1.
For $\al\in\AJhat$ we may assume, after possibly making a basis transformation,
that $e_1$ is a $J$-invariant vector in $\FSH^\al$.

\proclaim{Lemma 5.3}
Let $(\FSA,J)$ be a generalized quantum Gelfand pair.
Let $\al\in\Ahat$. The subspace of $J$-biinvariant elements in
\rom{Span}$\{t_{ij}^\al\}$ is equal to $\CC\, t_{11}^\al$ if $\al\in\AJhat$
and equal to $\{0\}$ otherwise.
\endproclaim

\demo{Proof}
Let $\al\in\Ahat$, $f\in\FSA^*$ and let the $c_{ij}$ be
arbitary complex coefficients. Then
$$
\eqalignno{
f.\Bigl(\sum_{i,j}c_{ij}\,t_{ij}^\al\Bigr)&=
\sum_{i,j,k}\pi_{kj}^\al(f)\,c_{ij}\,t_{ik}^\al,
\cr
\Bigl(\sum_{i,j}c_{ij}\,t_{ij}^\al\Bigr).f&=
\sum_{i,j,k}\pi_{ik}^\al(f)\,c_{ij}\,t_{kj}^\al.
\cr}
$$
Hence, $\sum_{i,j}c_{ij}\,t_{ij}^\al$ is $J$-biinvariant if and only if,
for all $f\in J$ and for all indices $k,l$ the following two equalities hold:
$$
\sum_j\pi_{kj}^\al(f)\,c_{lj}=0,\quad
\sum_i\pi_{ik}^\al(f)\,c_{il}=0.
$$
The first equality implies that $c_{lj}=0$ for all $l,j$ if $\al\notin\AJhat$
and that $c_{lj}=0$ for $l\ne1$ if $\al\in\AJhat$.
The second equality can be equivalently written as
$\sum_i\pi_{ki}^\al(f^*)\,\wb{c_{il}}=0$. Since $J=J^*$, we conclude that
also $c_{il}=0$ for $l\ne1$ if $\al\in\AJhat$.
Conversely, we see that $t_{11}^\al$ is $J$-biinvariant
if $\al\in\AJhat$.
\qed\enddemo

\proclaim{Proposition 5.4}
Let $(\FSA,J)$ be a generalized quantum Gelfand pair.
Then the statements of Proposition \rom{5.1} still hold if we read
$J$-biinvariant instead of $\FSB$-biinvariant and
$\AJhat$ instead of $\ABhat$.
\endproclaim

Again we call the elements $t_{11}^\al$ ($\al\in\AJhat$)
{\sl spherical elements} for the generalized quantum Gelfand pair
$(\FSA,J)$.
In order to obtain some analogue of Proposition 5.2, we need a
further assumption.

\proclaim{Proposition 5.5}
Let $(\FSA,J)$ be a generalized quantum Gelfand pair and assume that
$S(J)=J^-$.
\roster
\item"{(a)}"
If $a\in\FSA$ is $J$-biinvariant then so is $a^-$.
\item"{(b)}"
If $\al\in\AJhat$ then
$\wb\al\in\AJhat$ and $(t_{11}^\al)^-=t_{11}^{\wb\al}$ is a spherical element.
\endroster
\endproclaim

\demo{Proof}
Apply Lemma 3.4(c) with the involution $a\mapsto a^-$. We get:
$$
(f.\wb a)^-=(S(f))^-.a,\quad
(\wb a.f)^-=a.(S(f))^-\qquad(a\in\FSA,\;f\in J).
$$
Since $S(J)=J^-$, part (a) of the Proposition follows.
Next, part (b) is obtained because, for $\al\in\AJhat$,
$a:=(t_{11}^\al)^-$ is $\FSB$-biinvariant and satisfies $\ep(a)=1$
(cf.\ Proposition 5.4).
\qed\enddemo

\head 6. Discrete DJS-hypergroups from quantum Gelfand pairs \endhead
In this section it will be shown that, in case of a quantum Gelfand
pair of CQG algebras, one can associate
the structure of a discrete DJS-hypergroup with the linearization formula
(5.1). A similar result will be proved
in case of a generalized quantum Gelfand pair
$(\FSA,J)$ under the additional assumption that $S(J)=J^-$.

DJS-hypergroups were introduced by
Jewett \ccite{9}{pp.\ 12 and 17}, who called them convos.
Slightly different definitions were given almost simultaneously by 
Dunkl and Spector.
In many subsequent papers by various authors these structures were called
hypergroups. I follow the suggestion of G. Litvinov and K. Ross
to use the term DJS-hypergroups (after Dunkl, Jewett and Spector),
in order to distinguish these hypergroups from
the earlier introduced Delsarte-Levitan hypergroups.

Jewett's axioms for a DJS-hypergroup were neatly rephrased by 
Lasser \ccite{16}{\S1}, see his conditions (H1)--(H6). Below I will conform
to his terminology and notation.
Thus $K$ is a locally compact Hausdorff space, $M(K)$ denotes the space of all
complex regular Borel measures
and $M^1(K)$ the subset of all probability measures.
If $x\in K$ then $p_x$ denotes the corresponding point measure, i.e.,
$p_x\in M^1(K)$ and $p_x(\{x\})=1$.
A DJS-hypergroup is determined by $K$ together with the following three
data.
\roster
\item"{(H$*$)}"
A continuous mapping
$(x,y)\mapsto p_x*p_y\colon K\times K\to M^1(K)$ ({\sl convolution}),
where $M^1(K)$ bears the weak topology \wrt $C_c(K)$.
\item"{(H${}^-$)}"
An involutive homeomorphism $x\mapsto \wb x\colon K\to K$ ({\sl involution}).
\item"{(H$e$)}"
A fixed element $e\in K$ ({\sl unit element}).
\endroster

After identification of $x$ with $p_x$,
the mapping in (a) has a unique extension to a continuous bilinear mapping
$(\mu,\nu)\mapsto \mu*\nu\colon M(K)\times M(K)\to M(K)$.
The involution on $K$ gives rise to an involution $\mu\mapsto \mu^*$
on $M(K)$ defined by $\mu^*(E):=\wb{\mu(E^-)}$ ($E$ a Borel subset of $K$).

\definition{Definition 6.1}
Let the quadruple $(K,*,{}^-,e)$ be as above. Then this forms a
{\sl DJS-hypergroup} if the following conditions are satisfied.
\roster
\item"{(H1)}"
$p_x*(p_y*p_z)=(p_x*p_y)*p_z$ for all $x,y,z\in K$.
\item"{(H2)}"
$\supp(p_x*p_y)$ is compact for all $x,y\in K$.
\item"{(H3)}"
$(p_x*p_y)^-=p_{\,\wb y}*p_{\,\wb x}$ for all $x,y\in K$.
\item"{(H4)}"
$p_e*p_x=p_x=p_x*p_e$ for all $x\in K$.
\item"{(H5)}"
$e\in\supp(p_x*p_{\,\wb y})$ if and only if $x=y$.
\item"{(H6)}"
The mapping $(x,y)\mapsto\supp(p_x*p_y)$ of $K\times K$ into the space of
nonvoid compact subsets of $K$ is continuous, the latter space with the
topology as given in \ccite{9}{\S2.5}.
\endroster
\enddefinition

\proclaim{Theorem 6.2}
Let $(\FSA,\FSB)$ be a quantum Gelfand pair of CQG algebras.
Let $K:=\ABhat$, endowed with the discrete topology.
Put
$$
(p_\al*p_\be)(\ga):=c_{\al,\be}(\ga),\quad\al,\be,\ga\in K,
$$
where $c_{\al,\be}(\ga)$ is defined in \rom{(5.1)}.
Take $\al\mapsto\wb\al$ as the involutive mapping
defined just before Proposition \rom{5.2}.
Take for the element $e\in K$ the trivial \corep\ of $\FSA$.
Then the quadruple $(K,*,{}^-,e)$ forms a DJS-hypergroup.
\endproclaim

\demo{Proof}
\roster
\item"{(H$*$)}"
First we prove that $p_\al*p_\be\in M^1(K)$.
Indeed, $c_{\al,\be}(\ga)\ge0$ by Proposition 5.1(e),
and $(p_\al*p_\be)(K)=1$ because
$$
1=\ep(t_{11}^\al)\,\ep(t_{11}^\be)=
\sum_{\ga\in K}c_{\al\be}(\ga)\,\ep(t_{11}^\ga)=\sum_{\ga\in K} c_{\al\be}(\ga).
$$
The mapping $(\al,\be)\mapsto p_\al*p_\be$ is continuous since $K$ has the
discrete topology.
\item"{(H1)}"
Since $(t_{11}^\al\,t_{11}^\be)\,t_{11}^\ga=
t_{11}^\al\,(t_{11}^\be\,t_{11}^\ga)$, we have
$$
\sum_{\de\in K}c_{\be\ga}(\de)\,c_{\al\de}(\ze)=
\sum_{\de\in K}c_{\al\be}(\de)\,c_{\de\ga}(\ze),
$$
and this implies that $p_\al*(p_\be*p_\ga)=(p_\al*p_\be)*p_\ga$.
\item"{(H2)}"
The support of $p_\al*p_\be$ is compact since $c_{\al\be}(\ga)\ne0$ for only
finitely many~$\ga$.
\item"{(H3)}"
Since $(t_{11}^\al\,t_{11}^\be)^*=(t_{11}^\be)^*\,(t_{11}^\al)^*$,
we have that $\wb{c_{\al\be}(\wb\ga)}=c_{\,\wb\be\,\wb\al}(\ga)$,
and this implies
that $(p_\al*p_\be)^-=p_{\,\wb\be}*p_{\,\wb\al}$.
\item"{(H4)}"
Since $t_{11}^e\,t_{11}^\al=t_{11}^\al=t_{11}^\al\,t_{11}^e$,
we have that $c_{e\al}(\be)=\de_{\al\be}=c_{\al e}(\be)$,
and this implies that $p_e*p_\al=p_\al=p_\al*p_e$.
\item"{(H5)}"
Application of (2.17) to the identity
$t_{11}^\al\,t_{11}^{\wb\be}=\sum_\ga c_{\al\,\wb\be}(\ga)\,t_{11}^\ga$
yields
$$
h(t_{11}^\al\,t_{11}^{\wb\be}\,(t_{11}^\ga)^*)=
c_{\al\,\wb\be}(\ga)\,(F_\ga)_{11}/\tr F_\ga.
\tag6.1
$$
Hence $c_{\al\,\wb\be}(e)=h(t_{11}^\al\,t_{11}^{\wb\be})=
h(t_{11}^\al\,(t_{11}^\be)^*)$,
which is nonzero if and only if $\al=\be$. Now use that
$p_\al*p_{\,\wb\be}=\sum_\ga c_{\al\,\wb\be}(\ga)\,t_{11}^\ga$.
\item"{(H6)}"
The required continuity is immediate since $K$ has discrete topology.\qed
\endroster
\enddemo

\proclaim{Theorem 6.3}
Let $\FSA$ be a CQG algebra, $J\i\FSA^*$ a two-sided coideal satisfying
$J^*=J$ and $S(J)=J^-$, and
assume that $(\FSA,J)$ is a generalized quantum Gelfand pair.
Let $K:=\AJhat$, endowed with the discrete topology.
Take convolution, involution and unit element as in Theorem \rom{6.2}.
Then the quadruple $(K,*,{}^-,e)$ forms a DJS-hypergroup.
\endproclaim

\demo{Proof}
The proof is essentially the same as for Theorem 6.2, except for
property (H5).
After (6.1) we now continue with
$$
c_{\al\,\wb\be}(e)=h(t_{11}^\al\,t_{11}^{\wb\be})=
h(t_{11}^\al\,(t_{11}^\be)^-)>0,
$$
where we used Proposition 5.5 and (3.16).
\qed\enddemo

\head 7. The quantum group $SU_q(2)$ \endhead
By way of example of the theory developed in the previous sections we
treat now the CMQG algebra $\FSA=\FSA_q(SU(2))$ associated with the
quantum group $SU_q(2)$
(cf.\ for instance Woronowicz \cite{24}).
Fix $q\in(0,1)$.
Define $\FSA$ as the unital associative algebra with generators
$\al,\be,\ga,\de$ and relations
$$
\eqalignno{
\al\be&=q\be\al,\quad
\al\ga=q\ga\al,\quad
\be\de=q\de\be,\quad
\ga\de=q\de\ga,\quad
\be\ga=\ga\be,
\cr
\al\de-q\be\ga&=\de\al-q^{-1}\be\ga=1.
\cr}
$$
It turns out that $\FSA$ becomes a Hopf $*$-algebra under the following
actions of the comultiplication $\De\colon\FSA\ten\FSA$,
counit $\ep\colon\FSA\to\CC$
(unital multiplicative linear mappings),
antipode $S\colon\FSA\to\FSA$ (unital antimultiplicative linear mapping),
and involution $*\colon \FSA\to\FSA$ (unital antimultiplicative antilinear
mapping).
$$
\align
\De\albegade&=\albegade\ten\albegade,\quad
\ep\albegade=\left(\matrix1&0\\0&1\endmatrix\right),
\\
S\albegade&=\left(\matrix\de&-q^{-1}\be\\-q\ga&\al\endmatrix\right),\quad
\left(\matrix\al^*&\be^*\\ \ga^*&\de^*\endmatrix\right)=
\left(\matrix\de&-q\ga\\-q^{-1}\be&\al\endmatrix\right).
\endalign
$$
Here the formula for $\De$ has to be interpreted in the sense of
matrix multiplication: $\De(\al)=\al\ten\al+\be\ten\ga$, etc.

The matrix $u:=\albegade$ is an irreducible unitary matrix \corep\ of $\FSA$.
Hence $\FSA$ is a CMQG algebra.
Observe that
$$
u'':=S^2(u)=\left(\matrix\al&q^{-2}\be\\ q^2\ga&\de\endmatrix\right),
\text{ hence}\quad
F\,u=u''\,F,\text{ where}\quad
F:=\left(\matrix q^{-1}&0\\ 0&q\endmatrix\right).
$$
Since $\tr(F)=\tr(F^{-1})>0$, the operator $F$ is associated with $u$
as in Proposition 2.2.
Hence, we obtain by (3.9) that
$$
f_z(\al)=q^{-z},\quad
f_z(\de)=q^z,\quad
f_z(\be)=f_z(\ga)=0.
\tag7.1
$$
Then $f_z$ can be extended to $\FSA$ as a unital multiplicative
linear functional (cf.\ Proposition 3.3(b)).
For $f:=f_z$ the inner automorphism $\tau_f$ (cf.\ (3.3)),
acting on the generators, takes the form
$$
f_{-z}.\albegade.f_z=
\left(\matrix\al&q^{-2z}\be\\ q^{2z}\ga&\de\endmatrix\right).
$$
In particular (cf.\ Proposition 3.3(d)),
$$
f_{-1}.\albegade.f_1=\left(\matrix\al&q^{-2}\be\\ q^2\ga&\de\endmatrix\right)=
S^2\albegade.
$$
The second involution (3.12), acting on the generators, takes the
following explicit form.
$$
\left(\matrix\al^-&\be^-\\ \ga^-&\de^-\endmatrix\right)=
\left(\matrix\de&-\ga\\ -\be&\al\endmatrix\right).
$$

Next we summarize some results from Vaksman \& Soibelman \cite{22},
Masuda e.a.\ \cite{17}, and
Koornwinder \cite{10}.
Let the Hopf $*$-algebra $\FSB:=\FSA(U(1))$ be described as the algebra
generated by $z$ and $z^{-1}$ with relations $zz^{-1}=1=z^{-1}z$,
comultiplication $\De(z):=z\ten z$ and involution $z^*:=z^{-1}$.
Then $\FSB$ is a commutative CMQG algebra and $(\FSB,\Psi)$ is
a sub-CQG algebra of $\FSA=\FSA_q(SU(2))$, where
$$
\left(\matrix\Psi(\al)&\Psi(\be)\\ \Psi(\ga)&\Psi(\de)\endmatrix\right):=
\left(\matrix z&0\\ 0&z^{-1}\endmatrix\right).
$$

Up to equivalence, there is for each positive dimension $2l+1$
(where $l\in\Ahat=\{0,\thalf,1,\ldots\}$) a unique irreducible unitary matrix
\corep\ $t^l$ of $\FSA$.
Then $(t_{mn}^l)_{m,n=-l,-l+1,\ldots,l}$, the representative of $l\in\Ahat$,
can be chosen \st\
$\Psi(t_{mn}^l)=\de_{mn}\,z^{-2n}$.
We conclude that the pair $(\FSA,\FSB)$ is
a quantum Gelfand pair of CQG algebras
and that $\ABhat=\Zplus:=\{0,1,2,\ldots\}$.
Now $t_{00}^l$ is $\FSB$-biinvariant for $l\in\Zplus$ and
$$
t_{00}^l=p_l(\ga\ga^*;q^2),\quad\text{where}\quad
p_l(x;q):={}_2\phi_1(q^{-l},q^{l+1};q;q,qx),
\tag7.2
$$
a special case of the {\sl little $q$-Jacobi polynomials}
(cf.\ Andrews \& Askey \cite{2}).
Note that the $t_{00}^l$ mutually commute.
The involution $l\mapsto\wb l$ on $\ABhat$ is the identity mapping,
so $(t_{00}^l)^*=t_{00}^l$.

It was already pointed out in \ccite{12}{Example 7.7} that,
by application of (5.1),
the little $q$-Legendre polynomials have a linearization formula
$$
p_l\,p_m=\sum_k c_{lm}(k)\,p_k\quad
\text{with}\quad c_{lm}(k)\ge0.
\tag7.3
$$
It seems that this result has not yet been proved by analytic methods
(however, see the remark at the end of this section).
Theorem 6.2 gives the structure
of (commutative) DJS-hypergroup corresponding to (7.3).

Next we consider generalized quantum Gelfand pairs $(\FSA,J)$, where
we take $\FSA=\FSA_q(SU(2))$ and $J$ is a two-sided coideal in $\FSA^*$.
Results on this, obtained in Koornwinder \cite{11},
\cite{14}, \ccite{12}{\S9}, will be summarized and next an application
of Theorem 6.3 will be given.

We can completely characterize elements $A^z$ ($z\in\CC$), $B$ and $C$
in $\FSA^*$ by the properties
$$
\multline
A^z\albegade=
\left(\matrix q^{\half z}&0\\ 0&q^{-\half z}\endmatrix\right),\quad
B\albegade=\left(\matrix 0&1\\ 0&0\endmatrix\right),
\\
C\albegade=\left(\matrix 0&0\\ 1&0\endmatrix\right),
\endmultline
$$
and
$$
\De(A^z)=A^z\ten A^z,\quad
\De(B)=A\ten B+B\ten A^{-1},\quad
\De(C)=A\ten C+C\ten A^{-1}.
$$
Then
$$
\align
&A^zA^{z'}=A^{z+z'},\quad
A^0=\ep,\quad
A^zB=q^zBA^z,\quad
A^zC=q^{-z}CA^z,
\\
&BC-CB={A^2-A^{-2}\over q-q^{-1}}\,,\quad
\ep(A^z)=1,\quad
\ep(B)=0=\ep(C),
\\
&S(A^z)=A^{-z},\quad
S(B)=-q^{-1}B,\quad
S(C)=-qC,
\\
&(A^z)^*=A^{\wb z},\quad
B^*=C,\quad
C^*=B.
\endalign
$$
The elements $A,A^{-1},B,C$ generate the well-know quantized universal
enveloping algebra $\FSU_q(sl(2))$.
By (7.1) we have
$$
f_z=A^{-2z}.
$$
Hence, by (3.14),
$$
(A^z)^-=A^{\wb z},\quad
B^-=qC,\quad
C^-=q^{-1}B.
$$

For $\si\in\RR$ put
$$
X_\si:=iB-iC-{q^{-\si}-q^\si\over q^{-1}-q}\,(A-A^{-1}).
\tag7.4
$$
Then
$$
\De(X_\si)=A\ten X_\si+X_\si\ten A^{-1},\quad
(X_\si)^*=X_\si,\quad
S(X_\si)=-(X_\si)^-.
$$
Hence, $J:=\CC\,X_\si$ is a two-sided coideal in $\FSA^*$ satisfying
$J=J^*$ and $S(J)=J^-$.

The expression (7.4) for $X_\si$ was used in \ccite{12}{\S9}.
In \cite{14} a slightly different expression for $X_\si$ is used,
\st\ $(X_\si)^*=S(X_\si)$ instead of $(X_\si)^*=X_\si$.
However, it is easy to reformulate results from \cite{14} in terms of
(7.4).

It follows from \cite{14} that the pair $(\FSA,J)$ is a generalized
quantum Gelfand pair and that, with $\Ahat=\{0,\thalf,1,\ldots\}$ as above,
$\AJhat=\Zplus$.
The spherical element corresponding to $n\in\AJhat$ is a positive multiple
of
$$
p_n(\rho_\si;-q^{2\si+1},-q^{-2\si+1},q,q\mid q^2),
\tag7.5
$$
where
$$
\rho_\si:=
\thalf\bigl(\al^2+\be^2+\ga^2+\de^2+iq^\half\,(q^{-\si}-q^\si)\,
(\de\ga+\be\al-\de\be-\ga\al)
+(q^{-\si}-q^\si)^2\,\be\ga\bigr),
$$
and $p_n$ in (7.5) is an {\sl Askey-Wilson polynomial}
\cite{3} generally defined by
$$
p_n(\cos\th;a,b,c,d\mid q):=
a^{-n}(ab,ac,ad;q)_n\,
{}_4\phi_3\left[{
q^{-n},q^{n-1}abcd,ae^{i\th},ae^{-i\th}\atop ab,ac,ad};q,q\right].
$$
Now Propositions 5.4, 5.5 and Theorem 6.3 are applicable.
In particular, as earlier mentioned in \ccite{12}{\S9},
it follows from Proposition 5.4 that Askey-Wilson
polynomials with parameters as in (7.5) satisfy a linearization
formula (7.3) with nonnegative coefficients.
Until now no analytic proof has been published for this positivity result,
except for the case $\si=0$, which goes back to Rogers,
cf.\ Gasper \& Rahman \ccite{8}{\S8.5}.
However, it is quite possible that this result will follow by dualization
and analytic continuation of the explicit product formula for Racah
polynomials with two parameters equal, cf.\ Gasper \& Rahman
\ccite{7}{(1.14), (1.15), (1.10)}.
(I thank G. Gasper for this observation.)$\;$
Note also that the positivity of the coefficients in (7.3) for the case
of the little $q$-Legendre polynomials (7.2) is a limit case
of the analogous result for the Askey-Wilson polynomials with parameters
as in (7.5), see the limit formula in \ccite{14}{\S6}.

\head 8. Noumi's quantum analogue of the Gelfand pair $(U(N), SO(N))$ \endhead
A very interesting generalized quantum Gelfand pair quantizing
the pair $(U(N),\allowbreak SO(N))$ was recently studied by Noumi \cite{18}.
We will show that not just Proposition 5.4, but also
Proposition 5.5 and Theorem 6.3 are applicable to this
situation. We will heavily refer to \cite{18} for notation, see also
the summary in Floris \ccite{6}{\S2, Example 2}.

Fix $q\in(0,1)$.
The CMQG algebra $\FSA=\FSA_q(U(N))$ is generated by elements
$t_{ij}$ ($1\le i,j\le N$) satisfying relations \ccite{18}{(1.3)}
and by a central element $(\det_q)^{-1}$ which is the inverse of
$\det_q$ given by \ccite{18}{(1.5)}.
The Hopf $*$-structure is determined by requiring that both
$((\det_q)^{-1})$ and $t:=(t_{ij})_{i,j=1,\ldots,N}$
are unitary matrix \corep s of $\FSA$.
A Hopf subalgebra $\FSU=\FSU_q(gl(N))$ of $\FSA^*$
can be given with generators
$q^\la$ ($\la\in\ZZ^N$) and $e_k,f_k$ ($1\le k\le N-1$),
satisfying relations \ccite{18}{(1.8)} and with comultiplication, counit
and antipode given by \ccite{18}{(1.9)}.
Take $\ep_1,\ldots,\ep_N$ as a standard basis for $\ZZ^N$.
The generators of $\FSU$ are completely characterized as linear functionals
on $\FSA$
when it is known how the generators of $\FSA$ are evaluated by these linear
functionals. These evaluations are given  by \ccite{18}{(1.20)}
and can be reformulated as
$$
\align
&q^\la(t_{ij})=\de_{ij}\,q^{\lan\la,\ep_i\ran},\quad
e_i(t_{kl})=\de_{ik}\,\de_{i+1,l},\quad
f_i(t_{kl})=\de_{i+1,k}\,\de_{il},\tag8.1
\\
&q^{\ep_i}((\textstyle\det_q)^{-1})=q^{-1},\quad
e_i((\textstyle\det_q)^{-1})=f_i((\textstyle\det_q)^{-1})=0.
\endalign
$$

It follows from \ccite{19}{(3.1), (3.2)} that
$$
S^2(t_{ij})=q^{2i-2j}\,t_{ij}.
$$
Hence,
$$
F\,t=t''\,F,\quad\text{where}\quad
F_{ij}=\de_{ij}\,q^{-N-1+2i}.
$$
Since $\tr(F)=\tr(F^{-1})>0$, the operator $F$ is associated with $t$
as in Proposition 2.2.
Hence, we obtain by (3.9) that
$$
f_z(t_{ij})=\de_{ij}\,q^{(-N-1+2i)z},\quad
f_z((\textstyle\det_q)^{-1})=1,
$$
where the second formula follows because $(\det_q)^{-1}$ is group-like.
It follows that
$$
f_z=q^{2z\rho_0},\quad\text{where}\quad
\rho_0:=\sum_{i=1}^N q^{i-\half(N+1)}\,\ep_i,
\tag8.2
$$
and where we extend the first equality in (8.1) to arbitrary
$\la\in\CC^N$.
For such general $q^\la$ the relations \ccite{18}{(1.8)} still hold.
It follows by (3.14) and (8.2) that
$$
f^-=q^{\rho_0}\,f^*\,q^{-\rho_0},\qquad f\in\FSU.
\tag8.3
$$

Consider now the two-sided coideal $\gok_q(a)$ in $\FSU$ 
for the case (SO) with $a:=(1,1,\ldots,1)$, as defined in
\ccite{18}{(2.4)}. Then $(\gok_q(a))^*=\gok_q(a)$ (cf.\ \ccite{18}{(3.31)}).
Because of \ccite{18}{Prop.\ 2.4}, the smaller two-sided coideal $J$
spanned by the elements
$$
X_i:=q^{\ep_i}\,f_i-q^{\ep_{i+1}}\,e_i\quad(i=1,\ldots,N-1),
\tag8.4
$$
generates the same left and right ideal in $\FSU$ as $\gok_q(1,1,\ldots,1)$.
It follows from \ccite{18}{(1.28)} that
$$
X_i^*=-X_i.
\tag8.5
$$
Hence $J^*=J$. Now, by the results in \cite{18}, the pair $(\FSA,J)$
is a generalized quantum Gelfand pair.

Let us inspect if the condition $S(J)=J^-$ is satisfied.
On the one hand, by (8.4) and \ccite{18}{(1.9)},
$$
S(X_i)=-q^{-(\ep_i+\ep_{i+1})}\,(q\,q^{\ep_i}\,f_i-q^{-1}\,q^{\ep_{i+1}}\,e_i).
$$
On the other hand, by (8.5), (8.3) and \ccite{18}{(1.8)},
$$
(X_i)^-=-(q\,q^{\ep_i}\,f_i-q^{-1}\,q^{\ep_{i+1}}\,e_i).
$$
Hence $S(X_i)=q^{-(\ep_i+\ep_{i+1})}\,(X_i)^-$.
Since $q^{-(\ep_i+\ep_{i+1})}$ commutes with $e_i$ and $f_i$, we obtain
$$
(S(X_i))^-=q^{-(\ep_i+\ep_{i+1})}\,X_i=X_i\,q^{-(\ep_i+\ep_{i+1})}.
\tag8.6
$$
Inspection of the proof of Proposition 5.5 immediately yields
that the conclusions of this Proposition (and of Theorem 6.3)
still hold for the present pair $(\FSA,J)$.
Thus we can associate with it the structure of a discrete DJS-hypergroup.
By \cite{18} the hypergroup is commutative and the spherical elements
can be expressed in terms of Macdonald's symmetric $q$-polynomials for
root system $A_{N-1}$.

Noumi \cite{18} also considers a two-sided coideal in $\FSU$ which
quantizes the subgroup $Sp(n)$ (notated $Sp(2n)$ by Noumi)
of $U(2n)$ ($N=2n$).
Again, a generalized quantum Gelfand pair is obtained
(see also Floris \ccite{6}{\S2, Example 2}),
but I did not yet find a property similar to (8.6) for that case.

\remark{Note added in proof}
For CQG algebras see also, in addition to \cite{4} and \cite{15},
the preprint ``CQG algebras: a direct algebraic approach to compact
quantum groups'' by M. S. Dijkhuizen and T. H. Koornwinder
(Report AM-R9401, CWI, Amsterdam).

The positivity of the linearization coefficients in (7.3) for polynomials
as in (7.2) or (7.5) can also be obtained by expressing these
coefficients in terms of Clebsch-Gordan coefficients, see Vainerman
\cite{21} for the case (7.2) and the forthcoming preprint
``Askey-Wilson polynomials and the quantum $SU(2)$ group, survey and
applications'' by H. T. Koelink for the case (7.5).
\endremark

\enddocument